\input amstex
 \documentstyle{amsppt}
\magnification=1200
\hcorrection{.25in}
\advance\vsize-.75in

 \NoBlackBoxes

 \TagsOnRight

 \define\C{\Bbb C} 
 \define\Q{\Bbb Q} 
 \define\Z{\Bbb Z} 

 
 
 

\topmatter
 \title
On weights of $l$-adic representation 
 \endtitle
 \thanks Supported by the Russian Foundation for Fundamental Research (grant 97-01-00647). Research at MSRI is supported in part by NSF grant DMS-9022140. This paper was written at the period of my stay at the Institut des Hautes \' Etudes Scientifiques (Bures-sur-Yvette) and at the Mathematical Sciences Research Institute (Berkeley) in Spring of 1997. I would like to thank the members of both institutions for their warm hospitality and excellent working conditions. 
\endthanks
 \author S.G.  Tankeev
 \endauthor
 \address  Vladimir State University, ul. Gorkogo, 87,
Vladimir, 600026
, Russia
\endaddress
\email 
tankeev\%rtf\@vpti.vladimir.su
\endemail

 \abstract
Let $J$ be an abelian variety over a number field  
$k$, [$k:\Q$]$<\infty$.   Assume that 
${\text{Cent(End}}(J\otimes\overline k))=\Z$. If the division $\Q$-algebra ${\text{End}}^0(J\otimes\overline k)$ splits at a prime number $l$, then the $l$-adic representation is defined by the minuscule weights (microweights) of simple classical Lie algebras of types $A_m$, $B_m$, $C_m$ or $D_m$.
\endabstract

 \endtopmatter

 \document

\head
 Introduction
\endhead

0.1. Let $J$ be an abelian variety over a number field $k \subset \C ,
 [k:\Q] < \infty$. For a given prime number $l$ consider the natural $l$-adic representation 
$$
\rho _l: {\text {Gal}}(\overline{k}/k)\to\text{GL}(H^1_{et}(J\otimes\overline k,\Q_l)).
$$
We recall some conjectures and results concerning $\rho_l$.

0.2. {\bf Serre - Mumford - Tate conjecture }[5], [9]. {\it There exists  a canonical isomorphism of Lie algebras }
$$
{\text{Lie Im}}(\rho_l)\simeq{\text{Lie[MT}}(J\otimes_k\C)(\Q_l)],
$$
{\it where } ${\text{MT}}(J\otimes_k\C)$ {\it is the Mumford - Tate group associated to the canonical Hodge $\Q$-structure on $H^1(J\otimes_k\C,\Q)$.}

0.3. {\bf P.Deligne theorem }[14, Theorem 0.5.1b]. {\it Each simple factor $g$ of the reductive Lie algebra } ${\text{Lie MT}}(J\otimes_k\C)\otimes\C$ {\it is a classical Lie algebra of type } $A_m$, $B_m$, $C_m$ {\it or }$D_m$, {\it and the highest weight of any irreducible $g$-submodule }
\newline
$E\subset H^1(J\otimes_k\C,\Q)\otimes_{\Q}\C$ {\it is a minuscule weight }({\it microweight}). 

We recall the list of minuscule weights:

for type $A_m\,\,\,(m\geq 1)$ : $\omega_1,\,\omega_2,\,...,\,\omega_m$;

for type $B_m\,\,\,(m\geq 2)$ : $\omega_m$;

for type $C_m\,\,\,(m\geq 3)$ : $\omega_1$;

for type $D_m\,\,\,(m\geq 4)$ : $\omega_1,\,\omega_{m-1},\,\omega_m$;

for type $E_6$ : $\omega_1,\,\omega_6$;

for type $E_7$ : $\omega_7$,
\newline
where $\omega_1,\,\omega_2,\,...,\,\omega_m$ are fundamental weights [2, Chapter 8, \S 7, Section 3].

In virtue of the Deligne theorem the following conjecture may be considered as a weak form of the Serre - Mumford - Tate conjecture:

0.4. {\bf Conjecture on minuscule weights }[14, Conjecture 0.4]. {\it Let $J$ be an abelian variety over a number field $k\subset\C$. Then each  simple factor $g$ of the reductive Lie algebra }${\text{Lie Im}}(\rho_l)\otimes\overline\Q_l$ {\it is a classical Lie algebra of type $A_m$,\,$B_m$,\,$C_m$ or $D_m$, and the highest weight of any irreducible $g$-submodule $E\subset H^1_{et}(J\otimes\overline k,\Q_l)\otimes\overline\Q_l$ is a minuscule weight }({\it microweight}).

0.5. If ${\text{End}}(J\otimes\overline k)=\Z$ and ${\text{dim}}_k\,J$ is an odd integer, then by the well known J.-P.Serre theorem [3, Section 6.1] there exists a canonical isomorphism of the $l$-adic Lie algebras 
$$
{\text{Lie Im}}(\rho_l)\simeq {\text{Lie [MT}}(J\otimes_k\C)(\Q_l){\text{]}}.
$$ 
In particular, the Serre - Mumford - Tate conjecture and the conjecture on minuscule weights hold for $J$. The analogous results are obtained for some abelian varieties of another types [10], [11], [12]. 

In this paper we prove the following main theorem.

0.6. {\bf Theorem }. {\it Let $J$ be a simple abelian variety  over a number field $k$, \, 
\newline
$[k:\Q]<\infty$. Assume that }${\text{Cent(End}}(J\otimes\overline k))=\Z$. {\it If the division $\Q$-algebra }${\text{End}}^0(J\otimes\overline k)$ {\it splits at a prime number $l$, then each simple factor $g$ of the reductive Lie algebra  }${\text{Lie Im}}(\rho_l)\otimes\overline\Q_l$ {\it is  a simple Lie algebra of type $A_m$, \, $B_m$, \, $C_m$  \,or \,  $D_m$, and the highest weight of any irreducible $g$-submodule $E\subset H^1_{et}(J\otimes\overline k,\Q_l)\otimes\overline\Q_l$ is a minuscule weight }({\it microweight}).

I am grateful to J.-P.Serre, D.Bertrand and Y.Andr\' e for the stimulating discussions on this topic. 
 
\head
\S 1. Some remarks on $l$-adic representation
\endhead

1.1. In what lollows $V_l=H^1_{et}(J\otimes\overline k,\Q_l)$. Let $G_{V_l}$ denote the algebraic envelope of the compact $l$-adic Lie group ${\text{Im}}(\rho_l)$. It is well known that $G_{V_l}$ coincides with the closure of ${\text{Im}}(\rho_l)$ in the Zariski topology of the group ${\text{GL}}(V_l)$. 

 Fix a non-Archimedean place $v$ of the field $k$. Let $k(v)$ is its residue field of characteristic $p_v$, \,\,\,  $k_v$ is the completion of the field $k$ relative to the topology induced by the place $v$, \,\,\, $\overline v$ is any extension of $v$ to the field $\overline k$, \,\,\, $D(\overline v)$ is the decomposition group (an element $\sigma\in {\text{Gal}}(\overline k/k)$ belongs to $D(\overline v)\Leftrightarrow \sigma(\overline v)=\overline v\Leftrightarrow \forall \alpha\in\overline k\,\,\,\,\overline v(\sigma\alpha)=\overline v(\alpha)$), \,\,\, $I(\overline v)\subset{\text{Gal}}(\overline k/k)$ is the inertia subgroup (an element $\sigma\in {\text{Gal}}(\overline k/k)$ belongs to
 $I(\overline v)\Leftrightarrow \forall \alpha\in\overline k $ the relation $\overline v(\alpha)\geq 0 $ implies $\overline v(\sigma\alpha-\alpha)>0))$.

It is well known that $I(\overline v)$ is a normal subgroup in $D(\overline v)$,
and $D(\overline v)/I(\overline v)\simeq{\text{Gal}}(\overline{k(v)}/k(v))$. We may identify ${\text{Gal}}(\overline{k_v}/k_v)$ with the decomposition group $D(\overline v)$. 

1.2. Assume that $v$  is a non-Archimedean place of $k$ at which $J$ has a good reduction. Let $F_{\overline v}\in{\text{Gal}}(\overline k/k)$ be the Frobenius element. 

It is well known that the characteristic polynomial of $\rho_l(F_{\overline v}^{-1})$ coincides with the characteristic polynomial of the Frobenius endomorphism $\pi_v$ of the reduction $J_v$ of $J$ at $v$. We denote by $\Delta$ the set of all eigenvalues of $\rho_l(F_{\overline v}^{-1})$ (counting multiplicities). Let $\Gamma_v$ be a multiplicative subgroup of $\overline\Q^{\times}$ generated by $\Delta$.

By the Honda - Tate theory [13]  \, $\Q[\pi_v]=\prod K_i$, \, $K_i$ are number fields. The multiplicative group $\Q[\pi_v]^{\times}$ determines a $\Q$-torus $T_{\pi_v}=\prod R_{K_i/\Q}({\text{G}}_{mK_i})$, where $R_{K_i/\Q}$ are the Weil restrictions of scalar functors. Let $H_v$ be the smallest algebraic subgroup of $T_{\pi_v}$  defined over $\Q$,  such that $\pi_v\in H_v(\Q)$. As is well known, $H_v$  is a group of multiplicative type. The connected component of the identity in $H_v$ is called the {\it Frobenius torus } $T_v$. Moreover, $T_v$ is a $\Q$-torus  in $G_{V_l}$ [3, Section 3b].

1.3. We denote by $\C_{p_v}=(\overline\Q_{p_v})^{\wedge}$ the completion of the algebraic closure of the field $\Q_{p_v}$. According to Tate and Raynaud there exists the canonical Hodge - Tate decomposition 
$$
V_{\C_{p_v}}=V_{p_v}\otimes_{\Q_{p_v}}\C_{p_v}=V_{\C_{p_v}}(0)\oplus V_{\C_{p_v}}(-1),
$$
and dim\,$V_{\C_{p_v}}(0)$=dim\,$V_{\C_{p_v}}(-1)$  ([3, Section 4.3]).

 We define a morphism of algebraic $\C_{p_v}$-groups 
$$
h_{V_{p_v}}\,:\,{\text{G}}_{m\C_{p_v}}\to{\text{GL}}(V_{\C_{p_v}})
$$
by the formula 
$$
h_{V_{p_v}}(c)(x)=
\cases
x\,\,\,\,\,\,\,\,{\text{if}} \,\,\,\,\, x\in V_{\C_{p_v}}(0), \\
c^{-1}x\,\,\,\,\,\, {\text{if}} \,\,\,\,\, x\in V_{\C_{p_v}}(-1).
\endcases
$$

The algebraic envelope of $\rho_{p_v}(I(\overline v))$ is an algebraic subgroup of $G_{V_{p_v}}$. By the theorem of S.Sen [7], the 1-parameter subgroup $h_{V_{p_v}}$ is contained in the algebraic envelope of  $\rho_{p_v}(I(\overline v))$. In particular, one has ${\text{Im}}(h_{V_{p_v}})\subset G_{V_{p_v}}(\C_{p_v})$.

1.4. {\bf Lemma }. {\it Assume that } ${\text{Cent End}}^0(J\otimes\overline k)$ {\it is a totally real field, and $J(k)$ contains all points of order $p_v^2$ in $J(\overline k)$. Let $M\subset G_{V_{p_v}}$ be the smallest normal algebraic subgroup for which } ${\text{Im}}(h_{V_{p_v}})\subset M(\C_{p_v})$. {\it Then 
the canonical representation of the semisimple Lie algebra } ${\text{[Lie }}(M){\text{]}}^{ss}\otimes\overline\Q_{p_v}$ {\it in }  $V_{p_v}\otimes\overline\Q_{p_v}$ {\it is defined by the minuscule weights of the classical Lie algebras } ({\it of types } $A_m,\,B_m,\,C_m,\,D_m$).

{\it Proof }. By F.A.Bogomolov theorem [1]  Lie Im($\rho _{p_v}$)=Lie($G_{V_{p_v}}$) and
$G_{V_{p_v}}$ containes the group G$_m$ of homotheties. By G.Faltings
theorems [4] $G_{V_{p_v}}$ is reductive. In the case under consideration $G_{V_{p_v}}$ is connected [3, Proposition 3.6], and 
$$
G_{V_{p_v}}={\text{G}}_m\cdot S_{V_{p_v}},
$$
where $S_{V_{p_v}}$=[$G_{V_{p_v}},
G_{V_{p_v}}$] is the commutator subgroup of $G_{V_{p_v}}$ [3, Section 1.2.2b]. 
It is clear that  $M$ is a connected reductive algebraic group over $\Q_{p_v}$ with a nontrivial semisimple part $M^{ss}=[M,M]$ (because ${\text{Im }}(h_{V_{p_v}})$ is not contained in ${\text{G}}_m(\C_{p_v})$); the canonical representation $V_{p_v}$ of $M$ is faithful; if $N$ is a normal connected algebraic subgroup of $M$, defined over $\Q_{p_v}$, such that ${\text{Im }}(h_{V_{p_v}})\subset N(\C_{p_v})$, then $N$ is a normal algebraic subgroup of $G_{V_{p_v}}$ ( indeed, let $f:\,S_{V_{p_v}}^{sc}\to S_{V_{p_v}}$ be the universal covering morphism; since $M^{ss}$ is a normal algebraic subgroup of $S_{V_{p_v}}$, the group $(f^{-1}(M^{ss}))^0$ is a normal algebraic subgroup of $S_{V_{p_v}}^{sc}$, and so we get the decomposition 
$$
S_{V_{p_v}}^{sc}=(f^{-1}(M^{ss}))^0\times H'',
$$
 where $(f^{-1}(M^{ss}))^0$ and $H''$ are simply connected semisimple $\Q_{p_v}$-groups;
on the other hand, $(f^{-1}(N^{ss}))^0$ is a normal connected subgroup of $(f^{-1}(M^{ss}))^0$, so we get the decompositions 
$$
(f^{-1}(M^{ss}))^0=(f^{-1}(N^{ss}))^0\times H'
$$ 
and 
$$
S_{V_{p_v}}^{sc}=(f^{-1}(N^{ss}))^0\times H'\times H'';
$$ 
consequently $N^{ss}$ is a normal algebraic subgroup of $S_{V_{p_v}}$, and $N$ is a normal algebraic subgroup of $G_{V{p_v}}$), so $N=M$; the action of ${\text{G}}_{m\C_{p_v}}$ on $V_{\C_{p_v}}=V_{p_v}\otimes_{\Q_{p_v}}\C_{p_v}$ defined by $h_{V_{p_v}}$ is of weights $0$  and $-1$. In virtue of the well known results of J.-P.Serre [8, \S 3] the representation $V_{p_v}$ of $M$ is defined by minuscule weights of the classical Lie algebras. Lemma 1.4 is proved.

1.5. {\bf Corollary }. {\it Assume that }  ${\text{Cent End}}^0(J\otimes\overline k)$ {\it is a totally real field. Then for each prime number $l$ the reductive  Lie algebra }  ${\text{Lie Im}}(\rho_l)\otimes\overline\Q_l$ {\it has a simple factor $g$ such that the canonical representation of $g$ in   $V_l\otimes\overline\Q_l$  is defined by the minuscule weights of the classical Lie algebras.}

Indeed, we can choose a non-Archimedean place $v$ such that $p_v=l$. Replacing $k$ by some finite algebraic extension, we may assume that $J(k)$ containes all points of order $p_v^2$ in $J(\overline k)$.

\head
\S 2. On ${\text{Gal}}(\overline\Q/\Q)$-invariant subsets of $\Delta\cdot\Delta$ 
\endhead

2.1. {\bf Proposition }. {\it Let $J$ be an abelian variety over a number field $k$. Suppose that } ${\text{Cent End}}(J\otimes\overline k)=\Z$ {\it and the division $\Q$-algebra } ${\text{End}}^0(J\otimes\overline k)$ {\it splits at the prime number $l$.

After replacing $k$ by some finite extension we may assume that for abelian variety $J$ and for each non-Archimedean place $v$ of $k$ from some set $S_{non-Arch}^{generic}(J,k,l)$ of places of Dirichlet's density $1$ the following conditions hold:} 

1) {\it for a fixed integer } $n\geq 2$ {\it such that } $l^n>(2\cdot
{\text{dim}}_k J)^2$, {\it the } $l^n$-{\it torsion points of
$J(\overline k)$ are rational points over $k$, and $k$ containes all the }
$(l^n)^{th}$ {\it roots of unity;}

2) $p_v=$char$(k(v))>(2\cdot{\text{dim}}_k J)^2$;

3) Norm$_{k/\Q}(v)=p_v$; \, $p_v$ {\it is unramified in $k$};

4) $J$ {\it has a good reduction at $v$};

5) $\Gamma_v$ {\it is torsion-free,}  $G_{V_l}$ {\it is connected and }
$\rho _l(F_{\overline v}^{-1})\in T_v(\Q)$;

6) {\it the Frobenius torus } $T_v$ {\it is a maximal torus of } $G_{V_l}$
{\it and}
$$
{\text{rank}}(\Gamma_v)={\text{dim}}(T_v)={\text{rank}}(G_{V_l}).
$$

{\it Proof }. ${\text{Cent End}}(J\otimes\overline k)=\Z$ by the condition of Proposition 2.1. By Albert's classification [6, \S 21]  ${\text{End}}(J\otimes\overline k)=\Z$ or ${\text{End}}^0(J\otimes\overline k)$ is a quaternion division algebra over $\Q$ which splits at $l$. So 
$$
{\text{End}}(J\otimes\overline k)\otimes\Q_l=
\cases
\Q_l\,\,\,\,\,\,\,\,\,\,\,\,\,\,\,\,\,\,\,{\text{if}} \,\,\, J \,\,\, {\text{is of type I by Albert's classification}}, \\
{\text{M}}_2(\Q_l)\,\,\,\,\,\, {\text{if}} \,\,\, J \,\,\, {\text{is of type II or III by Albert's classification}} .
\endcases
$$

After replacing $k$ by some finite extension we may assume that for abelian variety $J$ and for $k$ the condition (2.1.1) holds. Then the Dirichlet density of the set $S_{non-Arch}^{generic}(J,k,l)$ of all non-Archimedean places of $k$ satisfying the conditions (2.1.2) - (2.1.6) is equal to 1  [3, Lemma (2.1), Proposition (3.6), Corollary (3.8)].  
 On the other hand,
$$
G_{V_l}={\text{G}}_m\cdot S_{V_l},
$$
where $S_{V_l}$=[$G_{V_l},
G_{V_l}$] is the commutator subgroup of $G_{V_l}$ [3, Section 1.2.2b].
Moreover, $V_l$ is a direct some of 1 or 2 copies of absolutely irreducible  $G_{V_l}$-module $U_l$, and $U_l\otimes\C_l$ has the Hodge - Tate decomposition with usual properties [3, Section 4.3]. Proposition 2.1 is proved.

From now on we assume that $v\in S^{generic}_{non-Arch}(J,k,l)$.

2.2. Let $f:\,S^{sc}_{V_l}\to S_{V_l}$ be the universal covering morphism. An isogeny $f$ extends to an isogeny $f_{{\text{G}}_m}:\,{\text{G}}_m\times S^{sc}_{V_l}\to {\text{G}}_m\cdot S_{V_l}=G_{V_l}$, which is defined by the formula $f_{{\text{G}}_m}(a\times s)=a\cdot f(s)$ for $a\in{\text{G}}_m,\,\,\,s\in S^{sc}_{V_l}$.

Consider the exact sequence 
$$
1\to {\text{Ker}}(f_{{\text{G}}_m})\to{\text{G}}_m\times S^{sc}_{V_l}\to {\text{G}}_m\cdot S_{V_l}=G_{V_l}\to 1.
$$
It determines the exact sequence of commutative algebraic $\Q_l$-groups
$$
1\to {\text{Ker}}(f_{{\text{G}}_m})\cap [f_{{\text{G}}_m}^{-1}(T_v)]^0\to [f_{{\text{G}}_m}^{-1}(T_v)]^0{\buildrel \varphi\over\longrightarrow} T_v\to 1,
\tag2.2.1
$$
where $[f_{{\text{G}}_m}^{-1}(T_v)]^0$
is a maximal $\Q_l$-subtorus of ${\text{G}}_m\times S^{sc}_{V_l}$, which coincides with the connected component of the unity of $f_{{\text{G}}_m}^{-1}(T_v)$, \,\,
$\varphi$ is an isogeny of $\Q_l$-tori, \, ${\text{Ker}}(f_{{\text{G}}_m})\cap [f_{{\text{G}}_m}^{-1}(T_v)]^0={\text{Ker}}(\varphi)$ is a finite group, which is isomorphic to some subgroup of ${\text{Cent}}(S^{sc}_{V_l})$. We rewrite the sequence (2.2.1) as 
$$
1\to {\text{Ker}}(\varphi)\to [f_{{\text{G}}_m}^{-1}(T_v)]^0{\buildrel \varphi\over\longrightarrow} T_v\to 1.
\tag2.2.2
$$

2.3. Recall that $S_{V_l}\otimes\overline\Q_l$ is a $\overline\Q_l$-almost simple algebraic group  $\Leftrightarrow S_{V_l}\otimes\overline\Q_l$ is semisimple and does not have infinite normal algebraic $\overline\Q_l$-subgroups different from $S_{V_l}\otimes\overline\Q_l$. 

If $S_{V_l}\otimes\overline\Q_l$ is a $\overline\Q_l$-almost simple algebraic group, then by Corollary 1.5 the $l$-adic representation $\rho_l$ is defined by the minuscule weights of the classical Lie algebras. 

2.4. Assume that  $S_{V_l}\otimes\overline\Q_l$ is not a $\overline\Q_l$-almost simple algebraic group.  Since $S^{sc}_{V_l}\otimes\overline\Q_l$ is a simply connected algebraic $\overline\Q_l$-group, it is a product of two proper semisimple $\overline\Q_l$-subgroups:
$$
S^{sc}_{V_l}\otimes\overline\Q_l=S_1\times S_2.
$$
Consider the canonical projections ${\text{pr}}_i:\,{\text{G}}_m\times S_1\times S_2\to S_i$ ($i=1,\,2$) and  
\newline
${\text{pr}}_0:\,{\text{G}}_m\times S_1\times S_2\to {\text{G}}_m$.
Since $[f_{{\text{G}}_m}^{-1}(T_v\otimes\overline\Q_l)]^0$ is a maximal $\overline\Q_l$-subtorus of ${\text{G}}_m\times S_1\times S_2$, we have:
$$
[f_{{\text{G}}_m}^{-1}(T_v\otimes\overline\Q_l)]^0=
$$
$$
{\text{pr}}_0([f_{{\text{G}}_m}^{-1}(T_v\otimes\overline\Q_l)]^0)\times{\text{pr}}_1([f_{{\text{G}}_m}^{-1}(T_v\otimes\overline\Q_l)]^0)\times{\text{pr}}_2([f_{{\text{G}}_m}^{-1}(T_v\otimes\overline\Q_l)]^0),
$$
and  ${\text{pr}}_i([f_{{\text{G}}_m}^{-1}(T_v\otimes\overline\Q_l)]^0)$ is a maximal $\overline\Q_l$-subtorus of $S_i\,\,\,(i=1,2)$, \,\,\,
\newline
${\text{pr}}_0([f_{{\text{G}}_m}^{-1}(T_v\otimes\overline\Q_l)]^0)={\text{G}}_m$. We get the decomposition 
$$
[f_{{\text{G}}_m}^{-1}(T_v\otimes\overline\Q_l)]^0={\text{G}}_m\times{\text{pr}}_1([f_{{\text{G}}_m}^{-1}(T_v\otimes\overline\Q_l)]^0)\times{\text{pr}}_2([f_{{\text{G}}_m}^{-1}(T_v\otimes\overline\Q_l)]^0).
\tag2.4.1
$$

We can find an element $y_v\in([f_{{\text{G}}_m}^{-1}(T_v\otimes\overline\Q_l)]^0)(\overline\Q_l)$ such that $\varphi(y_v)=\rho_l(F_{\overline v}^{-1})$ (recall that $\rho_l(F_{\overline v}^{-1})\in T_v(\Q)\subset  T_v(\overline\Q_l)$).

So the decomposition (2.4.1) gives the decomposition $y_v=\lambda\times y_{v1}\times y_{v2}$, 
where 
$$
\lambda\in{\text{G}}_m({\overline\Q_l})={\overline\Q_l}^{\times},
$$
$$
y_{vi}\in {\text{pr}}_i([f_{{\text{G}}_m}^{-1}(T_v\otimes\overline\Q_l)]^0)(\overline\Q_l)\subset S_i(\overline\Q_l). 
$$

On the other hand,
$$
U_l\otimes \overline\Q_l=\overline\Q_l\otimes W_1\otimes W_2,
$$
where $W_1$ is an irreducible $S_1$-module, $W_2$ is an
irreducible $S_2$-module. Let
$$
\rho _0:\,{\text{G}}_m\otimes\overline\Q_l\to {\text{GL}}(\overline\Q_l)=\overline\Q_l^{\times},
$$
$$
\rho _i:S_i\to {\text{GL}}(W_i)\,\,\,(i=1,2)
$$
are the corresponding representations. We have a commutative diagram
$$
\matrix
({\text{G}}_m\otimes\overline\Q_l)\times S_1\times S_2
&{\buildrel {\rho_0\otimes\rho_1\otimes\rho_2} \over
\longrightarrow} 
&{\text{GL}}(\overline\Q_l\otimes W_1\otimes W_2) \\
\downarrow f_{{\text{G}}_m} &\ &\vert\vert \\
G_{V_l}\otimes {\overline\Q_l}\,\,\,\,\,\,\,\,\,\,\,\, &\subset
&{{\text{GL}}(\overline\Q_l\otimes W_1\otimes W_2)}.
\endmatrix
\tag2.4.2
$$

Since $f_{{\text{G}}_m}(y_{vi})\in T_v(\overline\Q_l)\subset{\text{GL}}(U_l)(\overline\Q_l)$, we can diagonalize the linear operator $f_{{\text{G}}_m}(y_{vi})$ (over $\overline\Q_l$). Let $\Delta_i$ be the set of its eigenvalues. The commutative diagram (2.4.2) shows that the set $\Delta$ of all eigenvalues of the linear operator 
\newline
$f_{{\text{G}}_m}(y_v)=\varphi(y_v)=\rho_l(F_{\overline v}^{-1})$ is equal to 
$$
\Delta=\lambda\cdot\Delta_1\cdot\Delta_2
\tag2.4.3
$$
with the natural action of ${\text{Gal}}(\overline\Q/\Q)$ on $\Delta$ (recall that the characteristic polynomial of $\rho_l(F_{\overline v}^{-1})$ belongs to $\Z[t]$).

Let $\Gamma_{\Delta_i}$ be the multiplicative subgroup of $\overline\Q^{\times}$ generated by $\Delta_i$, \, 
\newline
$\Gamma_{\{\lambda\}}=\{\lambda^m\,\vert\,m\in\Z\}$. We denote by ${\text{X}}^{\ast}({\text{pr}}_i([\varphi^{-1}(T_v\otimes\overline\Q_l)]^0))$ the group of characters of $\overline\Q_l$-torus  ${\text{pr}}_i([\varphi^{-1}(T_v\otimes\overline\Q_l)]^0)$. It is evident that  
$$
\Gamma_v\subset\Gamma_{\{\lambda\}}\cdot\Gamma_{\Delta_1}\cdot\Gamma_{\Delta_2},
$$
so
$$
{\text{rank}}\,\Gamma_v\leq\sum {\text{rank X}}^{\ast}({\text{pr}}_i([\varphi^{-1}(T_v\otimes\overline\Q_l)]^0))({\text{pr}}_i(y_v))\leq
$$
$$
\sum {\text{rank X}}^{\ast}({\text{pr}}_i([\varphi^{-1}(T_v\otimes\overline\Q_l)]^0)={\text{rank(G}}_m\times S_1\times S_2)={\text{rank}}(G_{V_l}).
\tag2.4.4
$$
Condition (2.1.6) implies that 
$$
{\text{rank X}}^{\ast}({\text{pr}}_i([\varphi^{-1}(T_v\otimes\overline\Q_l)]^0))({\text{pr}}_i(y_v))={\text{rank X}}^{\ast}({\text{pr}}_i([\varphi^{-1}(T_v\otimes\overline\Q_l)]^0)).
\tag2.4.5
$$
Since there exists the canonical exact sequence
$$
{\text{X}}^{\ast}({\text{pr}}_i([\varphi^{-1}(T_v\otimes\overline\Q_l)]^0))\to{\text{X}}^{\ast}({\text{pr}}_i([\varphi^{-1}(T_v\otimes\overline\Q_L)]^0))({\text{pr}}_i(y_v))\to 0,
$$
and ${\text{X}}^{\ast}({\text{pr}}_i([\varphi^{-1}(T_v\otimes\overline\Q_l)]^0))$ is torsion-free, we deduce from (2.4.5) that 
$$
\Gamma_{\Delta_i}=
{\text{X}}^{\ast}({\text{pr}}_i([\varphi^{-1}(T_v\otimes\overline\Q_l)]^0))({\text{pr}}_i((y_v))
$$
is torsion-free. So (2.4.4) implies that 
$$
\Gamma_{\{\lambda\}}\cdot\Gamma_{\Delta_1}\cdot\Gamma_{\Delta_2}=
\Gamma_{\{\lambda\}}\times\Gamma_{\Delta_1}\times\Gamma_{\Delta_2}.
\tag2.4.6
$$

On the other hand, $W_i$ is a symplectic or orthogonal $S_i$-module, so $W_i^{\ast}=W_i$. Hence $\Delta_i^{-1}=\Delta_i$.

The Galois group ${\text{Gal}}(\overline\Q/\Q)$ acts in a natural way on $\Delta$ and on $\Delta\cdot\Delta$. For each element $\eta\in\Delta\cdot\Delta$ we define a map $T_{\eta}\,:\,\Delta\to\overline\Q^{\times}$ by the formula $T_{\eta}(\delta)=\eta\delta^{-1}$. It is easy to see that for each element $\sigma\in{\text{Gal}}(\overline\Q/\Q)$  we have (counting multiplicities)
$$
{\text{Card}}(T_{\eta}(\Delta)\cap\Delta)=
{\text{Card}}(T_{\sigma(\eta)}(\sigma(\Delta))\cap\sigma(\Delta))=
{\text{Card}}(T_{\sigma(\eta)}(\Delta)\cap\Delta),
$$
so  for any constant $c$ the set 
$$
\{\eta\in\Delta\cdot\Delta\,\vert\,{\text{Card}}(T_{\eta}(\Delta)\cap\Delta)=c\} 
$$
is ${\text{Gal}}(\overline\Q/\Q)$-invariant.  It is easy to see that 
$$
\eta=\lambda^2\,\Leftrightarrow\,{\text{Card}}(T_{\eta}(\Delta)\cap\Delta)={\text{Card}}(\Delta).
$$
Consequently $\lambda^2\in\Q$, and $\lambda^2$ is of absolute value $p_v$ (by the condition (2.1.3) and the well known A.Weil theorem 
[6, \S 21, Theorem 4]),  so we get the equality $\lambda^2=\pm p_v$.

Let $\rho$  be a complex conjugation on $\overline\Q$  induced by some embedding $\overline\Q\to\C$. According to (2.4.6)  each element $\delta\in\Delta$ has the unique decomposition $\delta=\lambda\cdot\delta_1\cdot\delta_2$, where $\delta_i\in\Delta_i$. We claim that 
$$
\rho(\lambda\delta_1\delta_2)=\lambda\delta_1^{-1}\delta_2^{-1}.
\tag2.4.7
$$
Indeed, $\rho(\lambda\delta_1\delta_2)=\lambda\delta_1'\delta_2'$ for some $\delta_i'\in\Delta_i$. Then by the  A.Weil theorem and by relation (2.1.3) we have $\lambda\delta_1\delta_2\rho(\lambda\delta_1\delta_2)=p_v$, so  
$$
\lambda^2\delta_1\delta_2\delta_1'\delta_2'=p_v,\,\,\, \lambda^4(\delta_1\delta_1'\delta_2\delta_2')^2=p_v^2, \,\,\, p_v^2(\delta_1\delta_1'\delta_2\delta_2')^2=p_v^2, \,\,\, (\delta_1\delta_1'\delta_2\delta_2')^2=1.
$$

Note that $\lambda^2=\lambda\delta_1\delta_2\cdot\lambda\delta_1^{-1}\delta_2^{-1}\in\Gamma_v$, \, so $\delta_1\delta_1'\delta_2\delta_2'=\lambda\delta_1\delta_2\cdot\lambda\delta_1'\delta_2'/\lambda^2\in\Gamma_v$. Since $\Gamma_v$ is torsion free (by the condition (2.1.5)) the relation $(\delta_1\delta_1'\delta_2\delta_2')^2=1$ implies $\delta_1'\delta_2'=(\delta_1\delta_2)^{-1}$. So the claim (2.4.7) is proved. 

The A.Weil theorem and (2.4.7) imply that 
$$
\lambda^2=\lambda\delta_1\delta_2\cdot\rho(\lambda\delta_1\delta_2)=p_v.
\tag2.4.8
$$

2.5. {\bf Theorem }. {\it Assume that the decomposition } $\Delta=\lambda\cdot\Delta_1\cdot\Delta_2$ {\it corresponds to the decomposition  $S^{sc}_{V_l}\otimes\overline\Q_l=S_1\times S_2$, where $S_1$ and $S_2$ are semisimple algebraic $\overline\Q_l$-groups. Then for each $i$ \,\,\,\,\,$\Delta_i\cdot\Delta_i$ does not contain a } ${\text{Gal}}(\overline\Q/\Q)$-{\it invariant subset $\neq\{1\}$}({\it without counting multiplicities}).

{\it Proof }. Assume that $\Delta_1\cdot\Delta_1$ containes a ${\text{Gal}}(\overline\Q/\Q)$-invariant subset $B\neq\{1\}$ (without counting multiplicities). It is evident that 
$$
{\text{Card}}(B)<{1\over 2}\cdot {\text{Card}}(\Delta)\cdot{\text{Card}}(\Delta)\leq{1\over 2}\cdot(2\cdot{\text{dim}}_k\,J)^2<{1\over 2}\cdot l^n.
\tag2.5.1
$$
Since $\Delta_2^{-1}=\Delta_2$, we have  $1\in\Delta_2\cdot\Delta_2$, so 
$$
\lambda^2\cdot B\subset\lambda^2\cdot\Delta_1\cdot\Delta_1\subset\lambda^2\cdot\Delta_1\cdot\Delta_1\cdot\Delta_2\cdot\Delta_2=\Delta\cdot\Delta.
$$
Consequently $\sum_{z\in\lambda^2\cdot B}z\in\Z$.

Assume that $p_v$ does not divide $\sum_{z\in \lambda^2\cdot B}z$. Then for each place $w$ of the field $\overline\Q$ lying over $p_v$ we have 
$$
w\left(\sum_{z\in\lambda^2\cdot B}z\right)=0.
$$
It follows that there exists
$z_w\in\lambda^2\cdot B$ such that
$w(z_w)=0$. Hence for each element $\delta_2\in\Delta_2$ 
$$
0=w(z_w)={1\over 2}\{w(z_w\delta_2^2)+w(z_w\delta_2^{-2})\}.
$$
Since both summands in the last brackets are nonnegative in virtue of
the relations $z_w\delta_2^{\pm 2}\in\Delta\cdot\Delta$,
we have the equalities
$$
w(z_w\delta_2^2)=w(z_w\delta_2^{-2})=0.
$$
So $w(\delta_2)=0$ for {\it all } $w\vert p_v$.
It follows that $\delta_2$ is a root of $1$ [14, Sublemma 3.4.0]. Since $\delta_2\in\Delta_2$ is an arbitrary element, we conclude that the multiplicative group $\Gamma_{\Delta_2}$ generated by $\Delta_2$ is a finite subgroup of $\overline\Q^{\times}$. We know that $\Gamma_{\Delta_2}$ is torsion-free. So $\Gamma_{\Delta_2}=\{1\}$, and 
$$
{\text{rank}}\,\,\Gamma_v\leq 1+{\text{rank}}\,\,\Gamma_{\Delta_1}= 1+{\text{rank X}}^{\ast}({\text{pr}}_1([\varphi^{-1}(T_v\otimes\overline\Q_l)]^0))({\text{pr}}_1(y_v))\leq 
$$
$$
1+{\text{rank X}}^{\ast}({\text{pr}}_1([\varphi^{-1}(T_v\otimes\overline\Q_l)]^0))=1+{\text{rank}}(S_1)< {\text{rank (G}}_m\times S_1\times S_2)={\text{rank}}(G_{V_l})
$$
contrary to the condition (2.1.6).

Consequently $p_v$ divides $\sum_{z\in\lambda^2\cdot B}z$. Since $\lambda^2= p_v$ by (2.4.8) and each $b\in B$ is of absolute value 1 in virtue of the A.Weil theorem [6, \S 21, Theorem 4], we get the relation 
$$
\sum_{z\in\lambda^2\cdot B}z=p_v\cdot b_v,
$$
where $b_v\in\Z$ is of absolute value $\leq{\text{Card}}(B)<{1\over 2}\cdot l^n$ in virtue of (2.5.1).

By the condition (2.1.3) $p_v$ is unramified in $k$, and 
 Norm$_{k/\Q}(v)$ $=p_v$. Consequently  $p_v$
splits completely in $k$, and we have
$$
k\otimes_{\Q}\Q_{p_v}\simeq\Q_{p_v}\times ...\times \Q_{p_v}.
\tag2.5.2
$$
By the condition (2.1.1) $k$ containes all the $(l^n)^{th}$
roots of unity. On the other hand it is well known that
$(\Q_{p_v}^{\times})_{tors}\simeq{\Z/(p_v-1)\Z}$. It follows from
(2.5.2) that $l^n\vert(p_v-1)$ and hence
$$
p_v=1\,({\text{ mod }}l^n).
\tag2.5.3
$$
By the condition (2.1.1) all the $l^n$-torsion points of
$J(\overline k)$ are rational over $k$. It follows that $\rho_l\vert_{{
\text{Gal}}(\overline k/k)}$ is trivial mod $l^n$ and hence for each $\delta\in\Delta$ we have the relation $\delta=1\,({\text{mod }}l^n)$. The same relation holds  for each element of $\Delta\cdot\Delta$. 

The known relation $\lambda^2\cdot B\subset\Delta\cdot\Delta$ implies 
$$
\sum_{z\in\lambda^2\cdot B}z={\text{Card}}(B)\,({\text{ mod }}l^n).
\tag2.5.4
$$
In virtue of (2.5.3) we deduce from (2.5.4) that $b_v={\text{Card}}(B)\,({\text{ mod }}l^n)$. It is clear that the inequalities  $|b_v|\leq {\text{Card}}(B)<{1\over 2}\cdot l^n$
imply the relation
$$
b_v={\text{Card}}(B).
$$
Consequently all elements of the set $\lambda^2\cdot B$ are equal to $p_v=\lambda^2$ contradicting to the condition $B\neq\{1\}$. Theorem 2.5 is proved.

2.6. In virtue of the condition (2.1.5) $\Gamma_v$ is torsion-free. So we may consider the canonical embeddings
$$
\Delta\subset\Gamma_v\subset\Gamma_v\otimes\Q.
$$
In virtue of Corollary 1.5 there exists a simple simply connected algebraic $\overline\Q_l$-group  $S_2$ of classical type such that the highest weight of each irreducible $S_2$-submodule of $H^1_{et}(J\otimes\overline k,\Q_l)\otimes\overline\Q_l$ is a minuscule weight. 

2.7. {\bf Lemma }. {\it The set $\Delta_2$ contains no three distinct points which are lying on one affine line of the space $\Gamma_{\Delta_2}\otimes\Q$.}

{\it Proof }. Assume that $S_2$ is a simple group of type  $A_m\,\,(m\geq 1)$. Then in N.Bourbaki's notations 
$$
\omega_i=\varepsilon_1+...+\varepsilon_i-{i\over{m+1}}(\varepsilon_1+...+\varepsilon_{m+1}),
$$
Weyl's group $W(R)$ is the group of all permutations of $\{\varepsilon_1,...,\varepsilon_{m+1}\}$ [2, Chapter 6, Section 4.7], dim $E(\omega_r)={{m+1}\choose r}$ [2, Chapter 8, Table 2], $E(\omega_r)$ is symplectic or orthogonal $\Leftrightarrow r={{m+1}\over 2}$ [2, Chapter 8, Table 1], and in this case 
$$
\omega_r={1\over 2}(\varepsilon_1+...+\varepsilon_{(m+1)/2}-\varepsilon_{(m+3)/2}-...-\varepsilon_{m+1}),
$$
$$
{\text{ch }} E(\omega_r)={\text{ch }} E(\omega_{(m+1)/2})=\sum_{a_i\in\{\pm 1\} \atop {a_1+...+a_{m+1}=0}} e^{a_1{\varepsilon_1/ 2}+...+a_{m+1}{\varepsilon_{m+1}/2}}=
$$
$$
=\sum_{a_i\in\{\pm 1\} \atop {a_1+...+a_m\in\{\pm 1\}}} e^{a_1{{(\varepsilon_1-\varepsilon_{m+1})}/ 2}+...+a_m{{(\varepsilon_m-\varepsilon_{m+1})}/2}}.
$$
Since the representation of $S_2$ in $W_2$ is symplectic or orthogonal, we have the following relation:
$$
\Delta_2=\{\alpha_1^{a_1}...\alpha_m^{a_m}\,\vert\,a_j\in\{\pm 1\},a_1+...+a_m\in \{\pm 1\}\},
$$ 
where $\alpha_1$, $\alpha_2$, ..., $\alpha_m$ are multiplicatively independent algebraic numbers. Hence elements of $\Delta_2$ are vertices of some cube in  $\Gamma_{\Delta_2}\otimes\Q$. It is evident that no three such vertices are lying on one affine line of $\Gamma_{\Delta_2}\otimes\Q$.

Assume that $S_2$ is a simple group of type  $B_m\,\,(m\geq 2)$. Then in N.Bourbaki's notations dim $E(\omega_m)=2^m$ [2, Chapter 8, Table 2],
$$
\omega_m={1\over 2}(\varepsilon_1+...+\varepsilon_m),
$$
$$
{\text{ch }} E(\omega_m)=\sum_{a_i\in\{\pm 1\}} e^{a_1{\varepsilon_1/ 2}+...+a_m{\varepsilon_m/2}}
$$
[2, Chapter 6, Section 4.5], so we may assume that
$$
\Delta_2=\{\alpha_1^{a_1}...\alpha_m^{a_m}\,\vert\,a_j\in\{\pm 1\}\},
$$
where $\alpha_1$, $\alpha_2$, ..., $\alpha_m$ are multiplicatively independent algebraic numbers. Hence elements of $\Delta_2$ are vertices of some cube in  $\Gamma_{\Delta_2}\otimes\Q$, and no three such vertices are lying on one affine line of $\Gamma_{\Delta_2}\otimes\Q$.

Assume that $S_2$ is a simple group of type  $C_m\,\,(m\geq 2)$.  Then dim $E(\omega_1)=2m$ [2, Chapter 8, Table 2], $\omega_1=\varepsilon_1$,
$$
{\text{ch }}E(\omega_1)=\sum_{a_i\in\{\pm 1\} \atop {i\in \{1,...,m\}}} e^{a_i{\varepsilon_i}}
$$
[2, Chapter 6, Section 4.6]. So we may assume that
$$
\Delta_2=\{\alpha_j^{a_j}\,\vert\,a_j\in\{\pm 1\},j\in\{1,...,m\}\},
$$
where $\alpha_1$, $\alpha_2$, ..., $\alpha_m$ are multiplicatively independent algebraic numbers. Hence elements of $\Delta_2$ are vertices of some cube in  $\Gamma_{\Delta_2}\otimes\Q$, and no three such vertices are lying on one affine line of $\Gamma_{\Delta_2}\otimes\Q$.

Assume that $S_2$ is a simple group of type  $D_m\,\,(m\geq 3)$.  Then dim $E(\omega_1)=2m$ [2, Chapter 8, Table 2], \, $\omega_1=\varepsilon_1$,
$$
{\text{ch }}E(\omega_1)=\sum_{a_i\in\{\pm 1\} \atop {i\in \{1,...,m\}}} e^{a_i{\varepsilon_i}},
$$
$$
{\text{dim}}\, E(\omega_{m-1})={\text{dim}} \,E(\omega_m)=2^{m-1},
$$
$$
{\text{ch }}E(\omega_{m-1})=\sum_{a_i\in\{\pm 1\} \atop {{\text{Card}}(\{ i\mid  a_i=-1 \})=1({\text{ mod 2 }})}} e^{a_1{\varepsilon_1/ 2}+...+a_m{\varepsilon_m/2}},
$$
$$
{\text{ch }}E(\omega_m)=\sum_{a_i\in\{\pm 1\} \atop {{\text{Card}}(\{ i\mid  a_i=-1 \})=0({\text{ mod 2 }})}} e^{a_1{\varepsilon_1/2}+...+a_m{\varepsilon_m/2}}
$$
[2, Chapter 6, Section 4.8].  So we may assume that
$$
\Delta_2=\{\alpha_j^{a_j}\,\vert\,a_j\in\{\pm 1\},j\in\{1,...,m\}\},
$$
or
$$
\Delta_2=\{\alpha_1^{a_1}...\alpha_m^{a_m}\,\vert\,a_j\in\{\pm 1\},{\text{Card}}(\{ j \vert a_j=-1 \})={\text{1(mod 2)}}\},
$$
or
$$
\Delta_2=\{\alpha_1^{a_1}...\alpha_m^{a_m}(a_j\in\{\pm 1\},{\text{Card}}(\{ j \mid a_j=-1 \})={\text{0(mod 2)}}\},
$$
where $\alpha_1$, $\alpha_2$, ..., $\alpha_m$ are multiplicatively independent algebraic numbers. Hence in any case elements of $\Delta_2$ are vertices of some cube in  $\Gamma_{\Delta_2}\otimes\Q$, and no three such vertices are lying on one affine line of $\Gamma_{\Delta_2}\otimes\Q$. Lemma 2.7 is proved.

2.8. {\bf Lemma }. {\it If $x$,\,$y$,\,$z\in\Delta$ are three distinct points lying on some affine line $L\subset\Gamma_v\otimes\Q$, then $x=\lambda\cdot\delta_1'\cdot\delta_2$, \, $y=\lambda\cdot\delta_1\cdot\delta_2$, \, $z=\lambda\cdot\delta_1''\cdot\delta_2$, where $\delta_1'$, $\delta_1$, $\delta_1''\in\Delta_1$, and $\delta_2\in\Delta_2$.}

{\it Proof }. According to (2.4.6) we may consider the canonical projection 
$$
{\text{pr}}_{\Delta_2}:\,\Gamma_v\otimes\Q=(\Gamma_{\{\lambda\}}\otimes\Q)\times(\Gamma_{\Delta_1}\otimes\Q)\times(\Gamma_{\Delta_2}\otimes\Q)\to\Gamma_{\Delta_2}\otimes\Q.
$$
If $x$, \, $y$, \, $z$ are distinct elements of $\Delta$ lying on an affine line $L$, then ${\text{pr}}_{\Delta_2}(x)$,  ${\text{pr}}_{\Delta_2}(y)$
,  ${\text{pr}}_{\Delta_2}(z)$ are lying on  ${\text{pr}}_{\Delta_2}(L)$. Note that if $y$ is situated between $x$ and $z$, then ${\text{pr}}_{\Delta_2}(y)$
is situated between ${\text{pr}}_{\Delta_2}(x)$ and  ${\text{pr}}_{\Delta_2}(z)$. If  ${\text{pr}}_{\Delta_2}(L)$ is a line, then ${\text{pr}}_{\Delta_2}(x)$, ${\text{pr}}_{\Delta_2}(y)$, ${\text{pr}}_{\Delta_2}(z)$ are {\it distinct } elements of ${\text{pr}}_{\Delta_2}(\Delta)=\Delta_2$ lying on an affine line ${\text{pr}}_{\Delta_2}(L)$ contrary to Lemma 2.7. Consequently ${\text{pr}}_{\Delta_2}(L)$ is a single point of $\Delta_2$. Lemma 2.8 is proved.

\head
\S 3. Roots, weights and Frobenius eigenvalues
\endhead

3.1. Now we start to prove Theorem 0.6. 

Let $R$ be the root system of  ${\text{Lie Im}}(\rho_l)\otimes\overline\Q_l$ relative to a Cartan subalgebra $h$. We denote by $(\,\,,\,\,):\,h^{\ast}\times h\to\overline\Q_l$ the tautological pairing. Let $R^V$ be the dual root system.
The sets $R$ and $R^V$ are invariant under the natural action of the Weyl group $W$. We define the lattice of weights $P(R)$ in $h^{\ast}$:
$$
P(R)=\{\lambda\in h^{\ast}\,\vert\,(\lambda,\alpha^V)\in\Z\,\,\,\forall \alpha\in R\}.
$$

3.2. Let $\Omega(V_l\otimes\overline\Q_l)$ be the collection of weights of $h\subset{\text{Lie Im}}(\rho_l)\otimes\overline\Q_l$. We have 
$$
\Omega(V_l\otimes\overline\Q_l)\subset P(R)\subset h^{\ast}.
$$
The collection $\Omega(V_l\otimes\overline\Q_l)$ is $W$-invariant and $R$-saturated [2, Chapter 8, \S 7, Propositions 3, 4]. This means that, for a weight $\lambda\in\Omega(V_l\otimes\overline\Q_l)$ and a root $\alpha\in R$  such that  $m=(\lambda,\alpha^V)>0$, the collection $\Omega(V_l\otimes\overline\Q_l)$ containes the arithmetic progression $\lambda,\,\lambda-\alpha,...,\lambda-m\alpha$ of length $m$.  If $m=(\lambda,\alpha^V)<0$, then $\Omega(V_l\otimes\overline\Q_l)$ contains the arithmetic progression $\lambda$, $\lambda+\alpha$,...,$\lambda+(-m)\alpha$ of length $(-m)$.

3.3.  Recall that a nonzero weight $\lambda\in P(R)$ is called a minuscule weight (microweight)$\Leftrightarrow$$(\lambda,\alpha^V)\in\{0,1,-1\}$ for all $\alpha\in R$ [2, Chapter 8, \S 7.3]. 

If the canonical $l$-adic representation $\rho_l$  is not defined by the minuscule weights, then $\Omega(V_l\otimes\overline\Q_l)$ containes the arithmetic progression of type $\lambda$, $\lambda-\alpha$, $\lambda-2\alpha$ or $\lambda$, $\lambda+\alpha$, $\lambda+2\alpha$. In this case there exist distinct elements $\delta\neq\delta'\neq\delta''\in\Delta$ such that $\delta^2=\delta'\cdot\delta''$ \, (indeed, if $\lambda$, $\lambda-\alpha$, $\lambda-2\alpha\in\Omega(V_l\otimes\overline\Q_l)$, then  $2\cdot(\lambda-\alpha)=\lambda+(\lambda-2\alpha))$.  

3.4. {\bf Remark }. It is evident that $\delta$ is the middle point of the segment $[\delta',\delta'']\subset\Gamma_v\otimes\Q$. Hence $\delta'$, $\delta$ and $\delta''$ are lying on one affine line $L\subset\Gamma_v\otimes\Q$.

3.5. Consider the set 
$$
\Delta^{line}=\{(\delta',\delta,\delta'')\,\vert\,\delta\neq\delta'\neq\delta''\in\Delta,\,\delta^2=\delta'\cdot\delta''\}.
$$
It is evident that $\Delta^{line}\subset\Delta\times\Delta\times\Delta$ is invariant under the action of the Galois group ${\text{Gal}}(\overline\Q/\Q)$  because for each $\sigma\in{\text{Gal}}(\overline\Q/\Q)$ we have $(\sigma(\delta))^2=\sigma(\delta^2)=\sigma(\delta')\cdot\sigma(\delta'')$. 

The last  relation shows that the set 
$$
B=\{\delta'/\delta,\,\,\delta''/\delta=(\delta'/\delta)^{-1}\,\vert\,(\delta',\delta,\delta'')\in\Delta^{line}\}\subset\Gamma_v
$$
is invariant under the natural action of ${\text{Gal}}(\overline\Q/\Q)$. 

3.6. Let $(\delta',\delta,\delta'')\in\Delta^{line}$. According to Lemma 2.8 we have the following decompositions:
$$
\delta'=\lambda\cdot\delta_1'\cdot\delta_2,\,\,\,
\delta=\lambda\cdot\delta_1\cdot\delta_2,\,\,\,
\delta''=\lambda\cdot\delta_1''\cdot\delta_2,
\tag3.6.1
$$
where $\delta_1',\,\delta_1,\,\delta_1''\in\Delta_1$, \, $\delta_2\in\Delta_2$.

Let $b=\delta'/\delta\in B$. Then $b=\delta_1'/\delta_1$, and we get from (3.6.1) the following decompositions:
$$
\delta'=\lambda\cdot b\cdot\delta_1\cdot\delta_2,\,\,\,
\delta=\lambda\cdot\delta_1\cdot\delta_2,\,\,\,
\delta''=\lambda\cdot b^{-1}\cdot\delta_1\cdot\delta_2.
$$ 

Since $b\cdot\delta_1=\delta_1'\in\Delta_1$ and $\delta_1^{-1}\in\Delta_1$, we have the relation $b\in\Delta_1\cdot\Delta_1$. So $B\subset\Delta_1\cdot\Delta_1$. On the other hand, $B\neq\{1\}$, so we have a nontrivial ${\text{Gal}}(\overline\Q/\Q)$-invariant subset $B\subset\Delta_1\cdot\Delta_1$ contrary to Theorem 2.5.

Hence $\rho_l$ is defined by the minuscule weights. A simple factor of type $E_6$ can't appear because $E(\omega_{1,6})$ is not a symplectic or orthogonal representation [2, Chapter 8, Table 1]. We exclude the case of a simple factor of type $E_7$ in the next paragraph.

\head
\S 4. Exclusion of type $E_7$
\endhead

4.1. In Bourbaki's notations [2, Chapter 6, \S 4] a simple Lie algebra of type $E_7$ has the following system of roots:
$$
R=\{\pm\varepsilon_i\pm\varepsilon_j\,\,\,(1\leq i<j\leq 6),\,\,\,\pm(\varepsilon_7-\varepsilon_8),
$$
$$
\pm{1\over 2}\cdot(\varepsilon_7-\varepsilon_8+\sum_{i=1}^6(-1)^{\nu(i)}\varepsilon_i)\,\,\,(\sum_{i=1}^6\nu(i)=1\,(\,{\text{mod 2))}}\}.
$$
Base of the system of roots is the following set:
$$
\{\alpha_1={1\over 2}(\varepsilon_1+\varepsilon_8)-{1\over 2}(\varepsilon_2+\varepsilon_3+\varepsilon_4+\varepsilon_5+\varepsilon_6+\varepsilon_7);
$$
$$
\alpha_2=\varepsilon_1+\varepsilon_2;\,\,\,\alpha_3=\varepsilon_2-\varepsilon_1;\,\,\,\alpha_4=\varepsilon_3-\varepsilon_2;
$$
$$
\alpha_5=\varepsilon_4-\varepsilon_3;\,\,\,\alpha_6=\varepsilon_5-\varepsilon_4;\,\,\,\alpha_7=\varepsilon_6-\varepsilon_5\}.
$$

For a given $\alpha\in R$ consider  the reflection  
$$
S_{\alpha}(x)=x-2\cdot{(x,\alpha)\over(\alpha,\alpha)}\cdot\alpha.
$$
Since $(\alpha,\alpha)=2$ for all $\alpha\in R$, we get the following formula:
$$
S_{\alpha}(x)=x-(x,\alpha)\cdot\alpha.
$$
These reflections $S_{\alpha}\,\,\,(\alpha\in R)$ generate the Weyl group $W(R)$ of order $2^{10}\cdot 3^4\cdot 5\cdot 7$.

4.2. We want to find the set of all weights of irreducible symplectic representation $E(\omega_7)$ of degree $56=2^3\cdot 7$ [2, Chapter 8, Tables 1, 2]. Since 
$$
\omega_7=\varepsilon_6+{1\over 2}(\varepsilon_8-\varepsilon_7)={1\over 2}(2\alpha_1+3\alpha_2+4\alpha_3+6\alpha_4+5\alpha_5+4\alpha_6+3\alpha_7)
$$ 
is a munuscule weight (microweight), the set of all weights of $E(\omega_7)$ coincides with $W(R)\cdot\omega_7$ [2, Chapter 8, \S 7, Proposition 6].

It is easy to find the matrix $A=(a_{ij})=(S_{\alpha_i}(\alpha_j))$:
$$
A=\left(
\matrix
-\alpha_1&\alpha_2&\alpha_1+\alpha_3&\alpha_4&\alpha_5&\alpha_6&\alpha_7\\
\alpha_1&-\alpha_2&\alpha_3&\alpha_2+\alpha_4&\alpha_5&\alpha_6&\alpha_7\\
\alpha_1+\alpha_3&\alpha_2&-\alpha_3&\alpha_3+\alpha_4&\alpha_5&\alpha_6&\alpha_7\\
\alpha_1&\alpha_2+\alpha_4&\alpha_3+\alpha_4&-\alpha_4&\alpha_4+\alpha_5&\alpha_6&\alpha_7\\
\alpha_1&\alpha_2&\alpha_3&\alpha_4+\alpha_5&-\alpha_5&\alpha_5+\alpha_6&\alpha_7\\
\alpha_1&\alpha_2&\alpha_3&\alpha_4&\alpha_5+\alpha_6&-\alpha_6&\alpha_6+\alpha_7\\
\alpha_1&\alpha_2&\alpha_3&\alpha_4&\alpha_5&\alpha_6+\alpha_7&-\alpha_7  
\endmatrix
\right).
$$
We consider a positive root 
$$
\varepsilon_3+\varepsilon_6=\alpha_2+\alpha_3+2\alpha_4+\alpha_5+\alpha_6+\alpha_7
$$
and the corresponding reflection $S_{\varepsilon_3+\varepsilon_6}$. It is easy to see that 

$S_{\varepsilon_3+\varepsilon_6}(\alpha_1)=\alpha_1+\alpha_2+\alpha_3+2\alpha_4+\alpha_5+\alpha_6+\alpha_7$,

$S_{\varepsilon_3+\varepsilon_6}(\alpha_2)=\alpha_2$,

$S_{\varepsilon_3+\varepsilon_6}(\alpha_3)=\alpha_3$,

$S_{\varepsilon_3+\varepsilon_6}(\alpha_4)=-\alpha_2-\alpha_3-\alpha_4-\alpha_5-\alpha_6-\alpha_7$,

$S_{\varepsilon_3+\varepsilon_6}(\alpha_5)=\alpha_2+\alpha_3+2\alpha_4+2\alpha_5+\alpha_6+\alpha_7$,

$S_{\varepsilon_3+\varepsilon_6}(\alpha_6)=\alpha_6$,

$S_{\varepsilon_3+\varepsilon_6}(\alpha_7)=-\alpha_2-\alpha_3-2\alpha_4-\alpha_5-\alpha_6$.
\newline
So we get the following chains:

$
2\cdot\omega_7=2\alpha_1+3\alpha_2+4\alpha_3+6\alpha_4+5\alpha_5+4\alpha_6+3\alpha_7{\buildrel S_{\varepsilon_3+\varepsilon_6}\over\longrightarrow}
$

$
2\alpha_1+\alpha_2+2\alpha_3+2\alpha_4+3\alpha_5+2\alpha_6+\alpha_7{\buildrel S_{\alpha_1}\over\longrightarrow}
$

$
\alpha_2+2\alpha_3+2\alpha_4+3\alpha_5+2\alpha_6+\alpha_7{\buildrel S_{\alpha_3}\over\longrightarrow}
$

$
\alpha_2+2\alpha_4+3\alpha_5+2\alpha_6+\alpha_7{\buildrel S_{\alpha_5}\over\longrightarrow}
$

$
\alpha_2+2\alpha_4+\alpha_5+2\alpha_6+\alpha_7{\buildrel S_{\alpha_4}\over\longrightarrow}
$

$
\alpha_2+\alpha_5+2\alpha_6+\alpha_7{\buildrel S_{\alpha_2}\over\longrightarrow}
$

$
-\alpha_2+\alpha_5+2\alpha_6+\alpha_7{\buildrel S_{\alpha_6}\over\longrightarrow}
$

$
-\alpha_2+\alpha_5+\alpha_7;
$

$
\qquad 2\alpha_1+\alpha_2+2\alpha_3+2\alpha_4+3\alpha_5+2\alpha_6+\alpha_7{\buildrel S_{\alpha_4}\over\longrightarrow}
$

$
\qquad 2\alpha_1+\alpha_2+2\alpha_3+4\alpha_4+3\alpha_5+2\alpha_6+\alpha_7{\buildrel S_{\alpha_2}\over\longrightarrow}
$

$
\qquad 2\alpha_1+3\alpha_2+2\alpha_3+4\alpha_4+3\alpha_5+2\alpha_6+\alpha_7{\buildrel S_{\alpha_3}\over\longrightarrow}
$

$
\qquad 2\alpha_1+3\alpha_2+4\alpha_3+4\alpha_4+3\alpha_5+2\alpha_6+\alpha_7{\buildrel S_{\alpha_4}\over\longrightarrow}
$

$
\qquad 2\alpha_1+3\alpha_2+4\alpha_3+6\alpha_4+3\alpha_5+2\alpha_6+\alpha_7{\buildrel S_{\alpha_5}\over\longrightarrow}
$

$
\qquad 2\alpha_1+3\alpha_2+4\alpha_3+6\alpha_4+5\alpha_5+2\alpha_6+\alpha_7{\buildrel S_{\alpha_6}\over\longrightarrow}
$

$
\qquad 2\alpha_1+3\alpha_2+4\alpha_3+6\alpha_4+5\alpha_5+4\alpha_6+\alpha_7;
$

$
2\alpha_1+\alpha_2+2\alpha_3+2\alpha_4+3\alpha_5+2\alpha_6+\alpha_7{\buildrel S_{\alpha_5}\over\longrightarrow}
$

$
2\alpha_1+\alpha_2+2\alpha_3+2\alpha_4+\alpha_5+2\alpha_6+\alpha_7{\buildrel S_{\alpha_6}\over\longrightarrow}
$

$
2\alpha_1+\alpha_2+2\alpha_3+2\alpha_4+\alpha_5+\alpha_7{\buildrel S_{\alpha_7}\over\longrightarrow}
$

$
2\alpha_1+\alpha_2+2\alpha_3+2\alpha_4+\alpha_5-\alpha_7;
$

$
\qquad \alpha_2+2\alpha_3+2\alpha_4+3\alpha_5+2\alpha_6+\alpha_7{\buildrel S_{\alpha_4}\over\longrightarrow}
$

$
\qquad \alpha_2+2\alpha_3+4\alpha_4+3\alpha_5+2\alpha_6+\alpha_7{\buildrel S_{\alpha_2}\over\longrightarrow}
$

$
\qquad 3\alpha_2+2\alpha_3+4\alpha_4+3\alpha_5+2\alpha_6+\alpha_7;
$

$
\alpha_2+2\alpha_3+2\alpha_4+3\alpha_5+2\alpha_6+\alpha_7{\buildrel S_{\alpha_5}\over\longrightarrow}
$

$
\alpha_2+2\alpha_3+2\alpha_4+\alpha_5+2\alpha_6+\alpha_7{\buildrel S_{\alpha_6}\over\longrightarrow}
$

$
\alpha_2+2\alpha_3+2\alpha_4+\alpha_5+\alpha_7{\buildrel S_{\alpha_3}\over\longrightarrow}
$

$
\alpha_2+2\alpha_4+\alpha_5+\alpha_7{\buildrel S_{\alpha_4}\over\longrightarrow}
$

$
\alpha_2+\alpha_5+\alpha_7{\buildrel S_{\alpha_5}\over\longrightarrow}
$

$
\alpha_2-\alpha_5+\alpha_7;
$

$
\qquad 2\alpha_1+\alpha_2+2\alpha_3+4\alpha_4+3\alpha_5+2\alpha_6+\alpha_7{\buildrel S_{\alpha_3}\over\longrightarrow}
$

$
\qquad 2\alpha_1+\alpha_2+4\alpha_3+4\alpha_4+3\alpha_5+2\alpha_6+\alpha_7;
$

$
\alpha_2+2\alpha_3+2\alpha_4+\alpha_5+\alpha_7{\buildrel S_{\alpha_7}\over\longrightarrow}
$

$
\alpha_2+2\alpha_3+2\alpha_4+\alpha_5-\alpha_7{\buildrel S_{\alpha_3}\over\longrightarrow}
$

$
\alpha_2+2\alpha_4+\alpha_5-\alpha_7;
$

$
\qquad \alpha_2+\alpha_5+\alpha_7{\buildrel S_{\alpha_7}\over\longrightarrow}
$

$
\qquad \alpha_2+\alpha_5-\alpha_7.
$

The continuation of this procedure leads to old weights or to opposite weights. Consequently the set $\Omega(E(\omega_7))$ of all weights of $E(\omega_7)$ is equal to ${1\over 2}\Omega\cup (-{1\over 2}\Omega)$, where $\Omega$ is the set of all elements in the chains above. As a result, we obtain the following list of elements of $\Omega$:

$
2\alpha_1+3\alpha_2+4\alpha_3+6\alpha_4+5\alpha_5+4\alpha_6+3\alpha_7;
$

$
2\alpha_1+3\alpha_2+4\alpha_3+6\alpha_4+5\alpha_5+4\alpha_6+\,\,\,\alpha_7;
$

$
2\alpha_1+3\alpha_2+4\alpha_3+6\alpha_4+5\alpha_5+2\alpha_6+\,\,\,\alpha_7;
$

$
2\alpha_1+3\alpha_2+4\alpha_3+6\alpha_4+3\alpha_5+2\alpha_6+\,\,\,\alpha_7;
$

$
2\alpha_1+3\alpha_2+4\alpha_3+4\alpha_4+3\alpha_5+2\alpha_6+\,\,\,\alpha_7;
$

$
2\alpha_1+3\alpha_2+2\alpha_3+4\alpha_4+3\alpha_5+2\alpha_6+\,\,\,\alpha_7;
$

$
2\alpha_1+\,\,\,\alpha_2+4\alpha_3+4\alpha_4+3\alpha_5+2\alpha_6+\,\,\,\alpha_7;
$

$
2\alpha_1+\,\,\,\alpha_2+2\alpha_3+4\alpha_4+3\alpha_5+2\alpha_6+\,\,\,\alpha_7;
$

$
2\alpha_1+\,\,\,\alpha_2+2\alpha_3+2\alpha_4+3\alpha_5+2\alpha_6+\,\,\,\alpha_7;
$

$
2\alpha_1+\,\,\,\alpha_2+2\alpha_3+2\alpha_4+\,\,\,\alpha_5+2\alpha_6+\,\,\,\alpha_7;
$

$
2\alpha_1+\,\,\,\alpha_2+2\alpha_3+2\alpha_4+\,\,\,\alpha_5\qquad \,\,\,\,\,+\,\,\,\alpha_7;
$

$
2\alpha_1+\,\,\,\alpha_2+2\alpha_3+2\alpha_4+\,\,\,\alpha_5\qquad \,\,\,\,\,-\,\,\,\alpha_7;
$

$
\qquad\,\,\,\,\,3\alpha_2+2\alpha_3+4\alpha_4+3\alpha_5+2\alpha_6+\,\,\,\alpha_7;
$

$
\qquad\,\,\,\,\,\,\,\,\alpha_2+2\alpha_3+4\alpha_4+3\alpha_5+2\alpha_6+\,\,\,\alpha_7;
$

$
\qquad\,\,\,\,\,\,\,\,\alpha_2+2\alpha_3+2\alpha_4+3\alpha_5+2\alpha_6+\,\,\,\alpha_7;
$

$
\qquad\,\,\,\,\,\,\,\,\alpha_2+2\alpha_3+2\alpha_4+\,\,\,\alpha_5+2\alpha_6+\,\,\,\alpha_7;
$

$
\qquad\,\,\,\,\,\,\,\,\alpha_2+2\alpha_3+2\alpha_4+\,\,\,\alpha_5\qquad\,\,\,\,\,+\,\,\,\alpha_7;
$

$
\qquad\,\,\,\,\,\,\,\,\alpha_2+2\alpha_3+2\alpha_4+\,\,\,\alpha_5\qquad\,\,\,\,\,-\,\,\,\alpha_7;
$

$
\qquad\,\,\,\,\,\,\,\,\alpha_2\qquad\,\,\,\,\,+2\alpha_4+3\alpha_5+2\alpha_6+\,\,\,\alpha_7;
$

$
\qquad\,\,\,\,\,\,\,\,\alpha_2\qquad\,\,\,\,\,+2\alpha_4+\,\,\,\alpha_5+2\alpha_6+\,\,\,\alpha_7;
$

$
\qquad\,\,\,\,\,\,\,\,\alpha_2\qquad\,\,\,\,\,+2\alpha_4+\,\,\,\alpha_5\qquad\,\,\,\,\,+\,\,\,\alpha_7;
$

$
\qquad\,\,\,\,\,\,\,\,\alpha_2\qquad\,\,\,\,\,+2\alpha_4+\,\,\,\alpha_5\qquad\,\,\,\,\,-\,\,\,\alpha_7;
$

$
\qquad\,\,\,\,\,\,\,\,\alpha_2\qquad\qquad\,\,\,\,\,\,\,\,\,\,+\,\,\,\alpha_5+2\alpha_6+\,\,\,\alpha_7;
$

$
\qquad\,\,\,\,\,\,\,\,\alpha_2\qquad\qquad\,\,\,\,\,\,\,\,\,\,+\,\,\,\alpha_5\qquad\,\,\,\,\,+\,\,\,\alpha_7;
$

$
\qquad\,\,\,\,\,\,\,\,\alpha_2\qquad\qquad\,\,\,\,\,\,\,\,\,\,+\,\,\,\alpha_5\qquad\,\,\,\,\,-\,\,\,\alpha_7;
$

$
\qquad\,\,\,\,\,\,\,\,\alpha_2\qquad\qquad\,\,\,\,\,\,\,\,\,\,-\,\,\,\alpha_5\qquad\,\,\,\,\,+\,\,\,\alpha_7;
$

$
\qquad\,\,\,\,\,\,\,\,\alpha_2\qquad\qquad\,\,\,\,\,\,\,\,\,\,-\,\,\,\alpha_5-2\alpha_6-\,\,\,\alpha_7;
$

$
\qquad\,\,\,\,\,\,\,\,\alpha_2\qquad\qquad\,\,\,\,\,\,\,\,\,\,-\,\,\,\alpha_5\qquad\,\,\,\,\,-\,\,\,\alpha_7.
$

4.3. Let $(x_1,...,x_7)$ be the coordinates on $\Q\alpha_1+...+\Q\alpha_7$ with respect to $\{\alpha_1,...,\alpha_7\}$. It is evident that the affine hyperplane $\{x\,\vert\,x_7=1\}$ contains 27 points of $\Omega\cup(-\Omega)$, and the affine hyperplane $\{x\,\vert\,x_7=-1\}$  containes 27 points of $\Omega\cup(-\Omega)$. Consequently there is a partition 
$$
\Omega\cup(-\Omega)=\{\pm(2\alpha_1+3\alpha_2+4\alpha_3+6\alpha_4+5\alpha_5+4\alpha_6+3\alpha_7)\}\cup
$$
$$
\{x\in\Omega\cup(-\Omega)\,\vert\,x_7=-1\}\cup\{x\in\Omega\cup(-\Omega)\,\vert\,x_7=1\}.
$$

Since $W(R)$ acts  transitively on $\Omega(E(\omega_7))={1\over 2}(\Omega\cup(-\Omega))$, we see that for each weight $\omega\in\Omega(E(\omega_7))$ there  exist parallel affine hyperplanes  $L_{-\omega}$, $L_{\omega}$ of dimension 6, which give a partition 
$$
\Omega(E(\omega_7))=\{\pm\omega\}\cup\{L_{-\omega}\cap\Omega(E(\omega_7))\}\cup\{L_{\omega}\cap\Omega(E(\omega_7))\},
$$
where 
$$
{\text{Card}}(L_{\pm\omega}\cap\Omega(E(\omega_7)))=27.
$$
These hyperplanes separate $-\omega$ and $\omega$ like on the following picture (where stars denote elements of $\Omega(E(\omega_7))$):
$$
\matrix
\ &\ &\ &\ &\ast&\vert &\ &\ &\ &\ast&\vert &\ &\ &\ \\
\ &\ &\ &\vert &\ast ...\ast &\vert &\ &\ &\vert &\ast ...\ast &\vert &\ &\ &\ \\
\ast &\ &\ &\vert &L_{-\omega} &\vert &\ &\ &\vert &L_{\omega} &\vert &\ &\ &\ast \\
-\omega &\ &\ &\vert &\ast ...\ast &\vert &\ &\ &\vert &\ast ...\ast &\vert &\ &\ &\omega \\
\ &\ &\ &\vert &\ast&\ &\ &\ &\vert &\ast&\ &\ &\ &\  
\endmatrix
\tag4.3.1
$$

4.4. {\bf Definition }. We say that a subset $D\subset\Delta$ has a configuration of type $E_7(\omega_7)\,\Leftrightarrow\,$ there exist a 7-dimensional affine subspace $L\subset\Gamma_v\otimes\Q$ and an affine isomorphism $f:\,L\simeq\Q\alpha_1+...+\Q\alpha_7$ such that the following conditions hold:

1) $D=L\cap\Delta$;

2) $f(D)=\Omega(E(\omega_7))$ is the system of all weights of the irreducible representation $E(\omega_7)$ of a simple Lie algebra of type $E_7$. 

4.5. {\bf Remark }. If $S_{V_l}^{sc}\otimes\overline\Q_l$ has a simple factor of type $E_7$, then there exists a subset $D\subset\Delta$ with a configuration of type $E_7(\omega_7)$. 

Indeed, if $\Delta_1'$ corresponds to a simple factor of type $E_7$, then due to 2.4  we get the decomposition 
$$
\Delta=\lambda\cdot\Delta_1'\cdot\Delta_1''\cdot\Delta_2.
$$ 
It is evident that for each $\delta_1''\in\Delta_1''$, $\delta_2\in\Delta_2$ the  affine space $L=(\Gamma_{\Delta_1'}\otimes\Q)\cdot\delta_1''\cdot\delta_2$ contains the subset $D=\Delta_1'\cdot\delta_1''\cdot\delta_2$ with a configuration of type $E_7(\omega_7)$.

4.6. {\bf Lemma }. {\it If $D\subset\Delta$ has a configuration of type $E_7(\omega_7)$, then } ${\text{pr}}_{\Delta_2}(D)$ {\it is a single point.}

{\it Proof }. We denote by $\Phi$ the composition of the following maps:
$$
\Q\alpha_1+...+\Q\alpha_7\,\,{\buildrel f^{-1}\over\longrightarrow}\,\,L\,\,{\buildrel{{\text{pr}}_{\Delta_2}}\over\longrightarrow}\,\,\Gamma_{\Delta_2}\otimes\Q.
$$
The map $\Phi$ is an affine map (possibly, degenerated). According to the results of Section 2.7 
$$
\Phi(\Omega(E(\omega_7)))={\text{pr}}_{\Delta_2}(f^{-1}(\Omega(E(\omega_7)))={\text{pr}}_{\Delta_2}(D)
$$
is a subset of the set of all vertices of a cube in $\Gamma_{\Delta_2}\otimes\Q$.

Assume that this subset contains at least two distinct points. Then we get from (4.3.1) the following picture (where affine subspaces $\Phi(L_{\pm\omega})$ are parallel and separate vertices $\Phi(-\omega)$ and $\Phi(\omega)$):
$$
\matrix
\ &\ &\ &\ &\ast &\vert &\ &\ &\ &\ast &\vert &\ &\ &\ \\
\ &\ &\ &\vert &\ast ...\ast &\vert &\ &\ &\vert &\ast ...\ast &\vert &\ &\ &\ \\
\ast &\ &\ &\vert &\Phi(L_{-\omega}) &\vert &\ &\ &\vert &\Phi(L_{\omega}) &\vert &\ &\ &\ast \\
\Phi(-\omega)&\ &\ &\vert &\ast ...\ast &\vert &\ &\ &\vert &\ast ...\ast &\vert &\ &\ &\Phi(\omega) \\
\ &\ &\ &\vert &\ast &\ &\ &\ &\vert &\ast &\ &\ &\ &\ 
\endmatrix
\tag4.6.1
$$
Since all stars are vertices of a cube, we can find an edge of a cube such that the following conditions hold: 1) $\Phi(-\omega)$ is lying on this edge; 2) the image of $\Phi(\omega)$ under the canonical projection of a cube onto this edge is not equal to $\Phi(-\omega)$. 

Let $\Phi'$ be the composition of $\Phi$ and the projection of a cube onto the  edge under consideration. We get from (4.6.1) the following picture:
$$
\matrix
\ast &\  &\ &\ &\ast &\  &\  &\ &\  &\ast &\  &\ &\ &\ast \\
\Phi'(-\omega)&\ &\ &\  &\Phi'(L_{-\omega})  &\  &\ &\ &\  &\Phi'(L_{\omega}) &\  &\ &\ &\Phi'(\omega) 
\endmatrix
$$
Consequently we arrive to a contradiction, because there are no ``intermediate'' vertices between two distinct vertices of a cube lying on one edge. Lemma 4.6 is proved.

4.7. {\bf Remark }. Assume that $D\subset\Delta$ has a configuration of type $E_7(\omega_7)$. Then for each element $\sigma\in{\text{Gal}}(\overline\Q/\Q)$ the set $\sigma(D)\subset\Delta$ has a configuration of type $E_7(\omega_7)$.

It is evident because ${\text{Gal}}(\overline\Q/\Q)$ acts on $\Gamma_v\otimes\Q$ as a group of affine transformations.

4.8. Consider the set 

$
\Delta^{E_7}=\{(\delta',\delta'',...,\delta^{(56)})\in\Delta\times ...\times\Delta=\Delta^{56}\,\vert \,\,\delta',\delta'',...\delta^{(56)} \,\,{\text{are distinct points }}
$

$
{\text{of }}\Delta  \,\,\,{\text{and the set }}\,\,\{\delta',\delta'',...,\delta^{(56)}\} \,\,\,{\text{has a configuration of type }}\,\,E_7(\omega_7)\}.
$

In virtue of Remark 4.7 this set is invariant under the action of the Galois group ${\text{Gal}}(\overline\Q/\Q)$. Consequenly the set 
$$
B=\{1,\delta''/\delta',...,\delta^{(56)}/\delta',\delta'/\delta'',...,\delta^{(56)}/\delta'',...\}\,\,\,\vert\,\,\,(\delta',...,\delta^{(56)})\in\Delta^{E_7}\}
$$
is invariant under the natural action of ${\text{Gal}}(\overline\Q/\Q)$. According to Lemma 4.6 we have the following decompositions:
$$
\delta'=\lambda\cdot\delta_1'\cdot\delta_2,\,\,\,
\delta''=\lambda\cdot\delta_1''\cdot\delta_2,\,\,\,
...\,\,\,
\delta^{(56)}=\lambda\cdot\delta_1^{(56)}\cdot\delta_2,
\tag4.8.1
$$
where $\delta_1',\,\,\delta_1'',\,\,...,\delta_1^{(56)}\in\Delta_1,\,\,\,\delta_2\in\Delta_2$.

Hence 
$$
B=\{1,\delta_1''/\delta_1',...,\delta_1^{(56)}/\delta_1',\delta_1'/\delta_1'',...,\delta_1^{(56)}/\delta_1'',...\}\,\,\,\vert(\delta',...,\delta^{(56)})\in\Delta^{E_7}\}
$$
is a ${\text{Gal}}(\overline\Q/\Q)$-invariant subset of $\Delta_1\cdot\Delta_1^{-1}=\Delta_1\cdot\Delta_1$, and $B\neq\{1\}$ contrary to Theorem 2.5. Thus a simple factor of type $E_7$ can't appear. Theorem 0.6 is proved.

\enddocument